\documentclass[10pt,twoside]{article}

\usepackage{epsfig}
\usepackage{latexsym}
\usepackage{amsfonts}
\usepackage{amsthm}

\DeclareSymbolFont{EulerScript}{U}{eus}{m}{n}
\DeclareSymbolFontAlphabet\mathscr{EulerScript}

\def\ps@myheadings{\let\@mkboth\@gobbletw
\def\@oddhead{\hbox{}
\rightmark\hfil\eightrm\thepage}   
\def\@oddfoot{}\def\@evenhead{\eightrm\thepage\hfil
\leftmark\hbox{}}\def\@evenfoot{}
\def\sectionmark##1{}\def\subsectionmark##1{}}
\def\runninghead#1#2{\pagestyle{myheadings}
\markboth{\centerline{\footnotesize\rm{#1}}}
{\centerline{\footnotesize\rm{#2}}}}

\def\abstract#1{{
	\centering{\begin{minipage}{4.2in}\footnotesize\baselineskip=10pt
	{\sc Abstract.\/} #1\par 
	\end{minipage}}\par}}

\renewenvironment{thebibliography}[1]
	{\frenchspacing
	\rm\baselineskip=11pt
	 \begin{list}{[\arabic{enumi}]}
	{\usecounter{enumi}\setlength{\parsep}{0pt}
	 \setlength{\leftmargin 13.7pt}{\rightmargin 0pt} %[1--9] ITEMS
	 \setlength{\itemsep}{2pt} \settowidth
	{\labelwidth}{[#1]}\sloppy}}{\end{list}}

\renewcommand{\section}[1] {\vspace{15pt}\addtocounter{section}{1}\setcounter{prop}{0}
\centerline{\bf\thesection. #1}\par\vspace{8pt}}
\newcommand{\nonumsection}[1] {\vspace{12pt}\centerline{\bf #1}
	\par\vspace{5pt}}

\newtheoremstyle{mythm}{10pt}{10pt}{\it}{}{\sc}{.}{ }{}
\theoremstyle{mythm}
\newtheorem{prop}{Proposition}[section]
\newtheorem{lemma}[prop]{Lemma}
\newtheorem{cor}[prop]{Corollary}
\newtheorem{thm}[prop]{Theorem}
\newtheorem*{Torresformula}{Torres formula}
\newtheorem*{Torres-Fox}{Fox-Torres relation}

\newcounter{itemlistc}
\newenvironment{itemlist}
    	{\setcounter{itemlistc}{0}
	 \begin{list}{-}
	{\usecounter{itemlistc}
	 \setlength{\parsep}{0pt}
	 \setlength{\itemsep}{5pt}}}{\end{list}}
	 
\newcommand{\Z}{{\mathbb Z}}
\newcommand{\lk}{{\ell k}}
\newcommand{\N}{{\mathscr N}}
\newcommand{\PP}{{\mathscr P}}
\newcommand{\E}{{\mathscr E}}
\newcommand{\m}{(\underline{m})}
\newcommand{\fb}{\overline{F}}
\newcommand{\inte}{\hbox{int}\,}
\newcommand{\rk}{\hbox{rk}\,}
\newcommand{\Ho}{\hbox{H}}
\newcommand{\egal}{\,\dot{=}\,}
\newcommand{\Hom}{\hbox{Hom}}

\begin{document}
	
\runninghead{THE MULTIVARIABLE ALEXANDER POLYNOMIAL}{DAVID CIMASONI}

\thispagestyle{empty}

\centerline{\bf STUDYING THE MULTIVARIABLE ALEXANDER}
\centerline{\bf POLYNOMIAL BY MEANS OF SEIFERT SURFACES}
\vskip1cm
\centerline{\sc david cimasoni \footnotetext[0]{{\em 2000 Mathematics Subject Classification:\/} 57M25.}
\footnotetext[0]{{\em Key words and phrases:\/} Alexander polynomial of a link, Seifert surface, Torres conditions.}}
\vskip1cm
\abstract{We show how Seifert surfaces, so useful for the understanding of the
Alexander polynomial $\Delta_L(t)$, can be generalized in order to study the multivariable Alexander polynomial
$\Delta_L(t_1,\dots,t_\mu)$. In particular, we give an elementary and geometric proof of the Torres formula.}
\vskip.5cm

\section{Introduction}
The technique of Seifert surfaces, discovered by {\sc Herbert Seifert} \cite{Sei} in 1935, enabled him to make great progress in the study
of the Alexander polynomial of a knot. In particular, he succeeded in characterizing among all Laurent polynomials $\Delta(t)$
those that can be realized as the Alexander polynomial of a knot.
The introduction by {\sc Ralph Fox} of the multivariable Alexander polynomial $\Delta_L(t_1,\dots,t_\mu)$ of a $\mu$-component oriented
link $L$ naturally gave rise to the corresponding question for this new invariant (see \cite[Problem 2]{Fox'}).
{\sc Guillermo Torres} made use of the free differential calculus -- developed at that time by Fox -- to give several conditions for a
polynomial $\Delta$ in $\Z[t^{\pm 1}_1,\dots,t^{\pm 1}_\mu]$ to be the Alexander polynomial of a $\mu$-component link \cite{Tor,F-T}.
Since then, very little progress has been made: it is known that the Torres conditions are not sufficient in general \cite{Hil,Pla},
but a complete algebraic characterization remains out of reach.

In this paper, we present an original approach to this problem. We show how the technique of Seifert surfaces can be generalized to
obtain a new geometric interpretation of $\Delta_L(t^{m_1},\dots,t^{m_\mu})$ for any integers $m_1,\dots,m_\mu$ (see Proposition
\ref{prop1}
and Corollary \ref{cor6}). If an equality holds for $\Delta_L(t^{m_1},\dots,t^{m_\mu})$ for any integers $m_1,\dots,m_\mu$, then it also
holds for $\Delta_L(t_1,\dots,t_\mu)$ (Lemma \ref{lemma2}); therefore, it is possible to prove properties of $\Delta_L$ with this method.
As an example, we give an elementary and geometric proof of the celebrated Torres formula,
valid for any link in a homology $3$-sphere. We also present several properties of $\Delta_L$ which
turn out to be equivalent to the Torres conditions (Proposition \ref{prop8}).
 
\section{Preliminaries}
Let us consider an oriented ordered link $L=L_1\cup\dots\cup L_\mu$ in a homology $3$-sphere $\Sigma$, and let $X$ be the exterior of $L$.
If $\widehat X\stackrel{\hat p}{\to}X$ denotes the universal abelian covering of $X$ and $\widehat X^0$ the inverse image by
$\hat p$ of a base point $X^0$ of $X$, the homology $\Ho_1(\widehat X,\widehat X^0)$ is endowed with a natural structure of a
module over the ring $\Lambda_\mu=\Z[t_1^{\pm 1},\dots,t_\mu^{\pm 1}]$. Given an $m\times n$ presentation matrix of
$\Ho_1(\widehat X,\widehat X^0)$ -- that is, the matrix $\PP$ corresponding to a presentation with $n$ generators and $m$
relations -- the $(n-i)\times(n-i)$ minor determinants of $\PP$ span an ideal of $\Lambda_\mu$ denoted by
$\E_i\Ho_1(\widehat X,\widehat X^0)$.
The greatest common divisor of these minor determinants is denoted by $\Delta_i\Ho_1(\widehat X,\widehat X^0)$; this
invariant is well defined up to multiplication by units of $\Lambda_\mu$, that is, by
$\pm t_1^{\nu_1}\cdots t_\mu^{\nu_\mu}$ with $\nu_i\in\Z$. In the sequel, we will write $\Delta\egal\Delta'$ if two elements 
$\Delta,\;\Delta'$ of a ring $R\,$ satisfy $\Delta=\varepsilon\Delta'$ for some unit $\varepsilon$ of $R$.
The Laurent polynomial $\Delta_1\Ho_1(\widehat X,\widehat X^0)$ is called the {\em Alexander polynomial\/} of the link $L\,$
\cite{Ale,Fox}. It is denoted by $\Delta_L(t_1,\dots,t_\mu)$.

Our method will be to prove statements on this polynomial in an indirect way, by studying all the
infinite cyclic coverings of $X$. Since these coverings are classified by
$\Hom(\Ho_1(X),\Z)\simeq\Ho^1(X;\Z)\simeq\Ho_1(L)=\bigoplus_{i=1}^\mu\Z L_i$, this
leads to the following definition \cite{E-N}. A {\em multilink\/} is an oriented link $L=L_1\cup\dots\cup L_\mu\,$ in a homology sphere $\Sigma$
together with an integer $m_i$ associated with each component $L_i$, with the convention that a component $L_i$ with multiplicity $m_i$ is the same
as $-L_i$ ($L_i$ with reversed orientation) with mutliplicity $-m_i$. Throughout this paper, we will write $\underline{m}$
for the ordered set of integers $m_1,\dots,m_\mu$, $d\,$ for their greatest common divisor, and $L\m$ for the multilink.
Finally, we will also denote by $\underline{m}$ the morphism $\Ho_1(X)\to\Z$ given by
$\underline{m}(\gamma)=\sum_{i=1}^\mu m_i\lk(L_i,\gamma)$. Let $\widetilde X\stackrel{\tilde p}{\to}X$
be the regular $\Z$-covering determined by $\underline{m}$. If $\widetilde X^0=\tilde p^{-1}(X^0)$,
the homology $\Ho_1(\widetilde X,\widetilde X^0)$ can be thought of as a module over the ring $\Z[t^{\pm 1}]$. The Laurent polynomial
$\Delta_{L\m}(t)=\Delta_1\Ho_1(\widetilde X,\widetilde X^0)$ is called the {\em Alexander polynomial\/} of the multilink $L\m$.
Note that if $\underline{m}\neq\underline{0}$, the exact sequence of the pair $(\widetilde X,\widetilde X^0)$ implies at once that
$\E_1\Ho_1(\widetilde X,\widetilde X^0)=\E_0\Ho_1(\widetilde X)$. Therefore, $\Delta_{L\m}(t)$ is also
equal to $\Delta_0\Ho_1(\widetilde X)$.

Here is the dictionary between the polynomials $\Delta_L$ and $\Delta_{L\m}$:

\begin{prop}[\rm Eisenbud-Neumann \cite{E-N}]\label{prop1}
$$
\Delta_{L\m}(t)\;\egal\cases{\Delta_L(t^{m_1}) & if $\mu=1$; \cr           
            	(t^d-1)\,\Delta_L(t^{m_1},\dots,t^{m_\mu}) & if $\mu\ge 2$.\cr}
$$ 
\end{prop}
\proof To check this equality, we need the well-known fact that $\E_1\Ho_1(\widehat X,\widehat X^0)=(\Delta_\ast)\cdot I$, where
$I\,$ is the augmentation ideal $(t_1-1,\dots,t_\mu-1)$ and $\Delta_\ast$ some polynomial in $\Lambda_\mu$. This can be proved by
purely homological algebraic methods using the fact that the group $\pi_1(X)$ has defect $\ge 1$ (see \cite[Theorem 6.1]{E-N}). By considering
a finite presentation of $\Ho_1(\widehat X,\widehat X^0)$ given by an equivariant cellular decomposition of $\widehat X$, it is easy to
show that $\Ho_1(\widehat X,\widehat X^0)\otimes_{\Lambda_\mu}\Z[t^{\pm 1}]=\Ho_1(\widetilde X,\widetilde X^0)$, where $\Z[t^{\pm 1}]$
is endowed with the structure of $\Lambda_\mu$-algebra given by $t_i\mapsto t^{m_i}$ for $i=1,\dots,\mu$. Hence,
\begin{eqnarray*}
\E_1\Ho_1(\widetilde X,\widetilde X^0)&=&\E_1(\Ho_1(\widehat X,\widehat X^0)\otimes_{\Lambda_\mu}\Z[t^{\pm 1}])\cr
	&=&(\Delta_\ast(t^{m_1},\dots,t^{m_\mu}))\cdot(t^{m_1}-1,\dots,t^{m_\mu}-1)\cr
	&=&(\Delta_\ast(t^{m_1},\dots,t^{m_\mu}))\cdot(t^d-1).
\end{eqnarray*}
Since $\Delta_L=(t_1-1)\,\Delta_\ast$ if $\mu=1$ and $\Delta_L=\Delta_\ast$ if $\mu\ge 2$, the proposition is proved. \endproof

In order to show that properties of $\Delta_{L\m}$ translate directly into properties of $\Delta_L$, we also need the following lemma.

\begin{lemma}\label{lemma2}
Consider two polynomials $\Delta$ and $\Delta'$ in $\Lambda_\mu$ such that
$$
\Delta(t^{m_1},\dots,t^{m_\mu})\egal\Delta'(t^{m_1},\dots,t^{m_\mu})\;\hbox{ in }\;\Z[t^{\pm 1}] 
$$
for all $(m_1,\dots,m_\mu)$ in $\Z^\mu$ except possibly a finite number of them. Then, $\Delta\egal\Delta'$ in $\Lambda_\mu$.
\end{lemma} 
\proof Without loss of generality, it may be assumed that $\Delta=\sum a_{i_1\cdots i_\mu} t_1^{i_1}\cdots t_\mu^{i_\mu}$ and
$\Delta'=\sum b_{j_1\cdots j_\mu} t_1^{j_1}\cdots t_\mu^{j_\mu}$ with $a_{0\dots 0}>0$, $b_{0\dots 0}>0$, and only non-negative indices
$i_k,j_k\ge 0$. By hypothesis, there are maps $\Z^\mu\stackrel{\varepsilon}{\to}\{\pm 1\}$ and $\Z^\mu\stackrel{\nu}{\to}\Z$
such that the equality
$$
\sum a_{i_1\cdots i_\mu} t^{\sum_km_ki_k}\;=\;\varepsilon(m_1,\dots,m_\mu)\,t^{\nu(m_1,\dots,m_\mu)}\,\sum b_{j_1\cdots j_n}t^{\sum_km_kj_k}
$$
holds for all but a finite number of $(m_1,\dots,m_\mu)$ in $\Z^\mu$. Let us choose an integer $N$ greater than $\max_k\deg_{t_k}\Delta$ and
$\max_k\deg_{t_k}\Delta'$, and set $m_1=1$, $m_2=N$, $\dots, m_\mu=N^{\mu-1}$. By choosing $N$ sufficiently large, it may be assumed
that the equality above holds for this ordered set of integers. Since all these integers are positive as well as
the coefficients $a_{0\dots 0}$ and $b_{0\dots 0}$, it follows that $\varepsilon(1,N,\dots,N^{\mu-1})=+1$ and
$\nu(1,N,\dots,N^{\mu-1})=0$. This gives
$$
\sum a_{i_1\cdots i_\mu} t^{\,i_1+Ni_2+\dots+N^{\mu-1}i_\mu}=\sum b_{j_1\cdots j_\mu} t^{\,j_1+Nj_2+\dots+N^{\mu-1}j_\mu}.
$$
But the equality $i_1+Ni_2+\dots+N^{\mu-1}i_\mu=j_1+Nj_2+\dots+N^{\mu-1}j_\mu$ with $0\le i_k,j_k<N\,$ for all $k\,$ implies that
$(i_1,\dots,i_\mu)=(j_1,\dots,j_\mu)$. Hence, $a_{i_1,\dots,i_\mu}=b_{i_1,\dots,i_\mu}$ for all multi-indices $(i_1,\dots,i_\mu)$, which
proves the result. \endproof

\section{Generalized Seifert surfaces}
One of the advantages of multilinks is that they can be studied via generalized Seifert surfaces \cite{E-N}. A {\em Seifert surface\/}
for a multilink $L\m$ is an open embedded oriented surface $F\subset\Sigma\setminus L$ such that, if $F_0$ denotes $F\cap(\Sigma\setminus\inte\N(L))$,
the closure $c\ell(F)$ of $F$ intersects a closed tubular neighborhood $\N(L_i)$ of $L_i$ as follows for each $i$:
\begin{itemlist}
\item{If $m_i \neq 0$, $c\ell(F)\cap\N(L_i)$ consists of $|m_i|$ sheets meeting along $L_i$; $F$ is oriented such that 
$\partial F_0=m_iL_i$ in $\Ho_1(\N(L_i))$.}
\item{If $m_i=0$, $c\ell(F)\cap\N(L_i)$ consists of discs transverse to $L_i$; $F$ is oriented such that the intersection number
of $L_i$ with each of these discs is the same (either always $+1$ or always $-1$).}
\end{itemlist}
\begin{figure}[Htb]
   \begin{center}
     \epsfig{figure=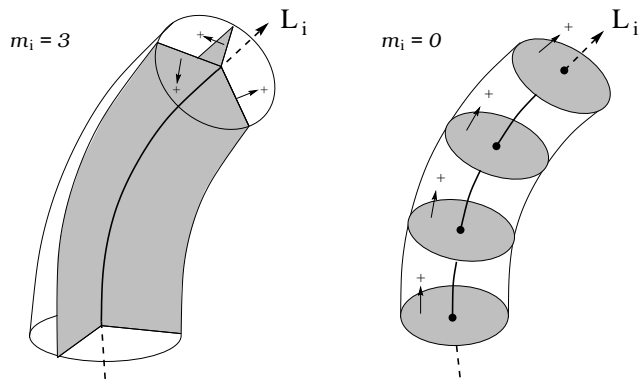,height=5cm}
     \caption{\footnotesize A Seifert surface near the multilink.}
     \label{fig1}
   \end{center}
\end{figure}
This is illustrated in Figure \ref{fig1}. 
Note that $F\subset\Sigma\setminus L$ and $F_0\subset\Sigma\setminus\inte\N(L)$ determine each other up to isotopy;
to simplify the notation, we will consider both of them as Seifert surfaces, and denote both by $F$. From now on, we will write $\fb$ for the union of
$F\subset\Sigma\setminus L$ and $L$.

\begin{lemma}[\rm Eisenbud-Neumann \cite{E-N}]\label{lemma3}
Let $F$ be a Seifert surface for a multilink $L\m$. Then, for $i=1\dots,\mu$, the intersection $F\cap\partial\N(L_i)$ gives
a $d_i$ component link which is the $(d_ip_i,d_iq_i)$-cable about $L_i$, where $p_i$ and $q_i$ are coprime, $d_ip_i=m_i$ and
$d_iq_i=-\sum_{j\neq i}m_j\lk(L_i,L_j)$.
\end{lemma}
\proof Let us denote by $(P_i,M_i)$ a basis of $\Ho_1(\partial\N(L_i))$ given by a standard parallel and meridian. Since $F\,$ is a Seifert
surface for $L\m$, $F\cap\partial\N(L_i)=m_iP_i+n_iM_i$ in $\Ho_1(\partial\N(L_i))$ for some integer $n_i$. Furthermore,
$\partial F=\sum_{j\neq i}m_jL_j+m_iP_i+n_iM_i$ in $\Ho_1(\Sigma\setminus\inte\N(L_i))$. By Alexander duality, this module is isomorphic
to $\Ho^1(\N(L_i);\Z)=\Z$, and the isomorphism is given by the linking number with $L_i$. It follows
$0=\lk(L_i,\partial F)=\sum_{j\neq i}m_j\lk(L_i,L_j)+n_i$, which gives the result. \endproof

In the usual case of an oriented link, a Seifert surface needs to be connected in order to be useful. In the general case of a multilink,
it has to be ``as connected as possible''. More precisely, a Seifert surface for $L\m$ is a {\em good Seifert surface\/} if it has
$d=\gcd\m$ connected components.

\begin{lemma}\label{lemma4}
Given a multilink $L\m$, there exists a good Seifert surface for $L\m$.
\end{lemma}
\proof One easily shows that there exists a Seifert surface for $L\m$ (see \cite[Lemma 3.1]{E-N}).
If $d>1$, a good Seifert surface for $L\m$ is given by $d\,$ parallel copies of a connected Seifert surface for $L(\frac{\underline{m}}{d})$.
Therefore, it may be assumed that $d=1$. Let $F\,$ be any Seifert surface for $L\m$ without closed component,
and let us denote by $i_+$ (resp. $i_-$) the epimorphism $\Ho_0(F)\to\Ho_0(\Sigma\setminus\fb)$ 
induced by the push in the positive (resp. negative) normal direction off $F$. If $i_+$ and $i_-$ are not isomorphisms, it is possible to
reduce the number of connected components of $F\,$ by handle attachement. So, let us assume that all the possible handle attachement(s)
have been performed, yielding $F=F_1\cup\dots\cup F_n$ with isomorphisms $i_+,i_-\colon\Ho_0(F)\to\Ho_0(\Sigma\setminus\fb)$. The
automorphism of $\Ho_0(F)$ given by $h=(i_-)^{-1}\circ i_+$ cyclically permutes the connected components of $F$. (Indeed, consider
a component $F_i$ of $F$; since $X=(\Sigma\setminus\fb)\cup F$ is path connected and $i_+,i_-$ are isomorphisms, there exists an
integer $m$ such that $F_i=h^m(F_1)$.) It easily follows that $\partial F_i=\partial F_j$ in $\Ho_1(\N(L))$ for $i,j=1,\dots,n$.
Therefore, the equality $\sum_{i=1}^\mu m_iL_i=\partial F=\sum_{j=1}^n\partial F_j=n\,\partial F_1$
holds in $\Ho_1(\N(L))=\bigoplus_{i=1}^\mu\Z L_i$.
Hence, $n$ divides $m_i$ for $i=1,\dots,\mu$. Since $\gcd(m_1,\dots,m_\mu)=1$, $F\,$ is connected.\endproof

Let us now turn to the natural generalization to multilinks of the Seifert form.
Given $F$ a good Seifert surface for $L\m$, the {\em Seifert forms\/} associated to $F$ are the bilinear forms
$$
\alpha_+,\alpha_-\/\colon\Ho_1(F)\times\Ho_1(\fb)\longrightarrow\Z
$$
given by $\alpha_+(x,y)=\lk(i_+x,y)$ and $\alpha_-(x,y)=\lk(i_-x,y)$,
where $i_+$ (resp. $i_-$) is the morphism $\Ho_1(F)\to\Ho_1(\Sigma\setminus\fb)$ induced by the push in the positive (resp.
negative) normal direction off $F$. (Note that we use the same notation for the morphisms $\Ho_0(i_\pm)$ and $\Ho_1(i_\pm)$; it will
always be clear from the context which dimension is concerned.) Let us denote by $A_+$ and $A_-$ matrices of these forms,
called {\em Seifert matrices}. Here is the generalization of Seifert's famous theorem.

\begin{thm}\label{thm}
Let $F\,$ be a good Seifert surface for $L\m$, and let $A_+,A_-$ be associated Seifert matrices. Then, $A_+-tA_-$ is a
presentation matrix of the module $\Ho_1(\widetilde X)$.
\end{thm}
\proof Given $F$ a good Seifert surface for $L\m$, let us denote $\Sigma\setminus\fb$ by $Y$. By the proof of Lemma \ref{lemma4}
it is possible to number the connected components $F=F_1\cup\dots\cup F_d\/$ and $Y=Y_1\cup\dots\cup Y_d\/$ such that
$i_+F_k=Y_k$ and $i_-F_k=Y_{k-1}$ (with the indices modulo $d\,$). Let us set $N=F\times(-1;1)$ an open bicollar of $F$,
$N_+=F\times(0;1)$, $N_-=F\times(-1;0)$ and $\{Y^i\}_{i\in \Z}$ (resp. $\{N^i\}_{i\in \Z}$) copies of $Y$ (resp. $N$).
Define
$$
E=\bigsqcup_{i\in\Z}{Y^i}\sqcup\bigsqcup_{i\in\Z}{N^i}{\Big /} \sim\,,
$$
where $Y^i\supset N_+\sim N_+\subset N^i\,$ and $Y^i\supset N_-\sim N_- \subset N^{i+1}$. The obvious projection
$E\stackrel{p}{\to}X$ is the infinite cyclic covering $\widetilde X\to X$ determined by $\underline{m}$.
Indeed, a loop $\gamma$ in $X$ lifts to a loop in $E$ if and only if the intersection number of $\gamma$ with $F$ is zero,
that is, if $\,0=\gamma\cdot F=\lk(L\m,\gamma)=\underline{m}(\gamma)$.

Consider the Mayer-Vietoris exact sequence of $\Z[t^{\pm 1}]$-modules associated to the decomposition
$\widetilde X={\big (}\bigcup_i Y^i{\big )} \cup {\big (}\bigcup_i N^i{\big )}$; it gives
\begin{eqnarray*}
(\Ho_1(F)\oplus \Ho_1(F))\otimes\Z[t^{\pm 1}]&\stackrel{\phi_1}{\to}&(\Ho_1(Y)\oplus\Ho_1(F))\otimes\Z[t^{\pm 1}]
\stackrel{\psi}{\longrightarrow}\Ho_1(\widetilde X)\to\\
(\Ho_0(F)\oplus\Ho_0(F))\otimes\Z[t^{\pm 1}]&\stackrel{\phi_0}{\to}&(\Ho_0(Y)\oplus\Ho_0(F))\otimes\Z[t^{\pm 1}],
\end{eqnarray*}
where the homomorphism $\phi_0$ is given by $(\alpha,\beta)\mapsto(i_+\alpha+t\,i_-\beta,\alpha+\beta)$.
Since $F$ is good, the homomorphisms $i_\pm\colon \Ho_0(F) \to \Ho_0(\Sigma \setminus\fb)$ are injective, and so is $\phi_0$.
Therefore, $\psi$ is surjective and
there is an exact sequence 
$$
(\Ho_1(F)\oplus\Ho_1(F))\otimes\Z[t^{\pm 1}]\stackrel{\phi_1}{\to}(\Ho_1(Y)\oplus \Ho_1(F))\otimes\Z[t^{\pm 1}]\to\Ho_1(\widetilde X)\to 0,
$$
with $\phi_1(\alpha,\beta)=(i_+\alpha+t\,i_-\beta,\alpha+\beta)$. This can be transformed into
$$
\Ho_1(F)\otimes\Z[t^{\pm 1}]\stackrel{\tilde\phi}{\to}\Ho_1(Y)\otimes\Z[t^{\pm 1}]\to\Ho_1(\widetilde X)\to 0,
$$
where $\tilde\phi(\alpha)=i_+\alpha-t\,i_-\alpha$. Let us fix basis $\cal B$ for $\Ho_1(F)$, $\overline{\cal B}$ for $\Ho_1(\fb)$, and
consider the basis $\overline{\cal B}^\ast$ for $\Ho_1(Y)$ which is dual to $\overline{\cal B}$ under Alexander duality. The matrix
of $i_+$ (resp. $i_-$) with respect to $\cal B$ and $\overline{\cal B}^\ast$ is given by $A^T_+$ (resp. $A^T_-$), where $A_+$ and $A_-$
are the Seifert matrices with respect to the basis $\cal B$ and $\overline{\cal B}$. Therefore, a matrix of $\tilde\phi$ is given by $A^T_+-tA^T_-$.
This concludes the proof. \endproof

\begin{cor}\label{cor6}
Let $L\m$ be a multilink with $\underline{m}\neq\underline{0}$.
If $m_i=\sum_{j\neq i}m_j\lk(L_i,L_j)=0$ for some index $i$, then $\Delta_{L\m}(t)=0$. If there is no such index, the matrices
$A_+$ and $A_-$ are square, and $\Delta_{L\m}(t)\egal\det(A_+-tA_-)$.
\end{cor}
\proof By the proof of Lemma \ref{lemma4}, a Seifert surface $F\,$ is good if and only if
$\rk\widetilde\Ho_0(F)=\rk\widetilde\Ho_0(\Sigma\setminus\fb)$ which is equal to $\rk\Ho_2(\fb)$ by Alexander duality. It is easy to show that
$\rk\widetilde\Ho_0(\fb)=r$, the number of indices $i$ with $m_i=\sum_{j\neq i}m_j\lk(L_i,L_j)=0$. Since $\chi(F)=\chi(\fb)$, it follows that
$\rk\Ho_1(\fb)=\rk\Ho_1(F)+r$. So if $r=0$, $A_+-tA_-$ is a square presentation matrix of $\Ho_1(\widetilde X)$ and if $r>0$,
it has more generators
than relations. It follows that $\Delta_0\Ho_1(\widetilde X)\egal\det(A_+-tA_-)$ if $r=0$, and $\Delta_0\Ho_1(\widetilde X)=0$ if
$r>0$. \endproof

\section{The Torres conditions}
Let us now illustrate how Corollary \ref{cor6}, along with Proposition \ref{prop1} and Lemma \ref{lemma2}, can be used to study the
multivariable Alexander polynomial. As an example, we present an elementary proof of the Torres formula \cite{Tor}, quite simpler
than the original proof. (On the other hand, it should be mentioned that more perspicuous proofs have since been given, for example
in \cite{Kaw}).

Throughout this section, we will denote by $\ell_{ij}$ the linking number $\lk(L_i,L_j)$.

\begin{figure}[Htb]
   \begin{center}
     \epsfig{figure=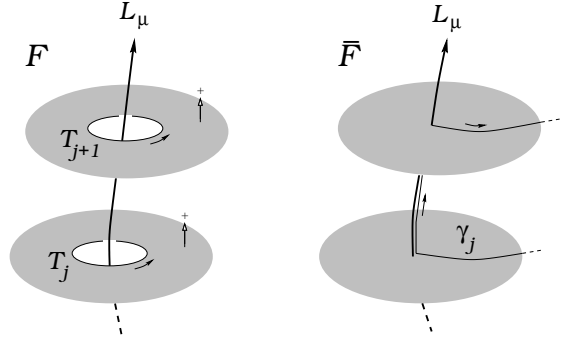,height=4.5cm}
     \caption{\footnotesize The proof of Lemma \ref{lemma7}.}
     \label{fig2}
   \end{center}
\end{figure}
\begin{lemma}\label{lemma7}
Let $L\m=L(m_1,\dots,m_{\mu-1},0)$ be a multilink, and let $L'(\underline{m}')=L'(m_1,\dots,m_{\mu-1})$ be the multilink obtained from
$L\m$ by removing the last component $L_\mu$. Then,
$$
\Delta_{L\m}(t)\;\egal\;(t^{\sum_im_i \ell_{i\mu}}-1)\,\Delta_{L'(\underline{m}')}(t).
$$
\end{lemma}
\proof If $m_i=\sum_{j\neq i}m_j\ell_{ij}=0$ for some index $i$, the lemma holds by Corollary \ref{cor6}. It may therefore be assumed
that there is no such index. Let $F$ be a good Seifert surface for $L\m$; then, a good Seifert surface for $L'(\underline{m}')$ is given by
$F'=F\cup(\fb\cap\N(L_\mu))$. By Lemma \ref{lemma3}, $\fb\cap\N(L_\mu)$ consists of $d_\mu=\sum_{i=1}^{\mu-1}m_i \ell_{i\mu}$ discs
(recall Figure \ref{fig1}). Furthermore, $\fb=\overline{F}'\cup L_\mu$. Therefore, we have the 
natural isomorphisms
$$
\Ho_1(F)=\Ho_1(F')\oplus\bigoplus_{j=1}^{d_\mu}\Z T_j\quad\hbox{and}\quad\Ho_1(\fb)=\Ho_1(\fb')\oplus\bigoplus_{j=1}^{d_\mu}\Z\gamma_j,
$$
where the cycles $T_j$ correspond to the boundaries of the discs, and the $\gamma_j$ are the transverse cycles depicted in Figure
\ref{fig2}. The associated Seifert matrices $A_\pm$ and $A'_\pm$ are related by
$$
A_+=\pmatrix{A'_+ & 0 \cr
 			\ast &\matrix{1& & \cr &\ddots&\cr & &1\cr}\cr}\quad \hbox{and}\quad
A_-=\pmatrix{A'_- & 0 \cr
 			\ast &\matrix{0& & &1\cr 1&0& &\cr &\ddots&\ddots&\cr & &1&0\cr}\cr}.
$$
Corollary \ref{cor6} then gives
\begin{eqnarray*}
\Delta_{L\m}&\dot{=}&\det(A_+-tA_-)\;=\;\left|\begin{array}{cc}
					A'_+-tA'_- &0\\
 					\ast &\matrix{1& & &-t\cr -t&1& & \cr &\ddots&\ddots& \cr & &-t&1\cr}
 					\end{array}\right|\\
 	&=&(t^{d_\mu}-1)\,\det(A'_+-tA'_-)\,\egal\,(t^{d_\mu}-1)\,\Delta_{L'(\underline{m}')}(t)
\end{eqnarray*}
and the lemma is proved.\endproof

The demonstration of the Torres formula is now a mere translation of Lemma \ref{lemma7} via Proposition \ref{prop1}.

\begin{Torresformula}{\rm \cite{Tor}}
Let $L=L_1\cup\dots\cup L_\mu$ be an oriented link with $\mu\ge 2$ components, and let $L'$ be the sublink $L_1\cup\dots\cup L_{\mu-1}$.
Then,
$$
\Delta_L(t_1,\dots,t_{\mu-1},1)\;\egal\cases{{{t_1^{\ell_{12}}-1}\over{t_1-1}}\,\Delta_{L'}(t_1) & if $\mu=2$;\cr           
            	(t_1^{\ell_{1\mu}}\cdots t_{\mu-1}^{\ell_{\mu-1,\mu}}-1)\,\Delta_{L'}(t_1,\dots,t_{\mu-1})& if $\mu>2$.\cr}
$$
\end{Torresformula}
\proof Let us denote by $\Delta'$ the right-hand side of this formula, and let $m_1,\dots,m_{\mu-1}$ be arbitrary integers with
$d=\gcd(m_1,\dots,m_{\mu-1})>0$. We have the equalities
\begin{eqnarray*}
\Delta'(t^{m_1},\dots,t^{m_{\mu-1}})&=&\cases{{{t^{m_1\ell_{12}}-1}\over{t^{m_1}-1}}\,\Delta_{L'}(t^{m_1}) & if $\mu=2$;\cr           
            	 (t^{\sum_i m_i\ell_{i\mu}}-1)\,\Delta_{L'}(t^{m_1},\dots,t^{m_{\mu-1}})&if $\mu>2$,\cr}\\
\hbox{(Proposition \ref{prop1})}&\dot{=}& {{1}\over{t^d-1}}\,(t^{\sum_i m_i\ell_{i\mu}}-1)\,\Delta_{L'(\underline{m}')}(t) \\
\hbox{(Lemma \ref{lemma7})}&\dot{=}& {{1}\over{t^d-1}}\,\Delta_{L\m}(t)\\
\hbox{(Proposition \ref{prop1})}&\dot{=}&\Delta_L(t^{m_1},\dots,t^{m_{\mu-1}},1)
\end{eqnarray*}
and the proof is settled by Lemma \ref{lemma2}. \endproof

\noindent Using the same method, it is not hard to show the following result.

\begin{Torres-Fox}{\rm \cite{Tor,F-T}}
Let $L=L_1\cup\dots\cup L_\mu$ be an oriented link with $\mu\ge 2$ components.
Then,
$$
\Delta_L(t_1^{-1},\dots,t_\mu^{-1})=(-1)^\mu\,t_1^{\nu_1-1}\cdots t_n^{\nu_\mu-1}\,\Delta_L(t_1,\dots,t_\mu)
$$ 
with integers $\nu_i$ such that $\nu_i\equiv \sum_{j}\ell_{ij}\!\!\!\pmod{2}$ if $\Delta_L\neq 0$.
\end{Torres-Fox}

These results provide necessary conditions for a polynomial $\Delta$ in $\Lambda_\mu$ to be the Alexander polynomial of a $\mu$-component link with fixed
$\lk(L_i,L_j)=\ell_{ij}$. They are known as the {\em Torres conditions} (see \cite{Lev} for a precise statement). Since these conditions are not sufficient
\cite{Hil,Pla}, the problem is now to find stronger conditions.
By means of a close study of the homology $\Ho_1(F)$ and $\Ho_1(\fb)$, it is possible to find necessary conditions for a polynomial
$\Delta$ in $\Z[t^{\pm 1}]$ to be the Alexander polynomial of a multilink. Via Proposition \ref{prop1}, this translates into the
following result (see \cite{cim} for a proof).
 
\begin{prop}\label{prop8}
Let $L\,$ be an oriented link with $\mu\ge 2$ components. Then, its Alexander polynomial $\Delta_L$ satisfies the following conditions.
For all integers $\underline{m}=(m_1,\dots,m_\mu)$ with $d=gcd(m_1,\dots,m_\mu)$ and $d_i=gcd(m_i,\sum_j m_j \ell_{ij})$, there exists
some polynomial $\nabla_{L\m}(t)$ in $\Z[t^{\pm d}]$ such that:
\begin{itemlist}
\item{$\prod_{i=1}^\mu(t^{d_i}-1)\,\nabla_{L\m}(t)\egal(t^d-1)^2\,\Delta_L(t^{m_1},\dots,t^{m_\mu});$}
\item{$\nabla_{L\m}(t^{-1})=\nabla_{L\m}(t)$;}
\item{$|\nabla_{L\m}(1)|={d^2\;D\over d_1\cdots d_\mu m_1\cdots m_\mu}$, where $D$ is any $(\mu-1)\times (\mu-1)$ minor
determinant of the matrix 
$$\pmatrix{
-\sum_j m_1m_j\ell_{1j} & m_1m_2\ell_{12} & \dots  & m_1m_\mu\ell_{1\mu} \cr
m_1m_2\ell_{12} & -\sum_j m_2m_j\ell_{2j} & \dots  & m_2m_\mu\ell_{2\mu} \cr
\vdots & \vdots  & \ddots & \vdots \cr
m_1m_\mu\ell_{1\mu} & m_2m_\mu\ell_{2\mu} & \dots & -\sum_j m_\mu m_j\ell_{\mu j}};$$}
\item{If $m_i=0$ for some index $i$, then $\nabla_{L\m}=\nabla_{L'(\underline{m}')}$, where $L'$ denotes the sublink $L\setminus L_i$ and
$\underline{m}'=(m_1,\dots,\widehat{m_i},\dots,m_\mu)$. \endproof}
\end{itemlist}
\end{prop}

This result easily implies the Torres conditions. It can also be thought of as a generalization of a theorem of Hosokawa \cite{Hos}, which
corresponds to the case $m_1=\dots=m_\mu=1$. At first sight, it might therefore seem more general than the Torres conditions.
Unfortunately, this is not the case: it can be shown that every polynomial $\Delta$ which
satisfies the Torres conditions also satisfies the conditions of Proposition \ref{prop8} (see \cite{cim}).

\smallskip

By means of a somewhat closer study of the Seifert matrices $A_\pm$, it should be possible to find new properties of $\Delta_{L\m}$.
They would translate into properties of $\Delta_L$, and provide new conditions, stronger than the ones of Torres.

\nonumsection{Acknowledgments}
I wish to express my thanks to Jerry Levine, to my advisor Claude Weber, and to Mathieu Baillif.

\vskip0.5cm
\noindent
\parbox{0.5\textwidth}{\footnotesize \sc David Cimasoni\\
Section de Math\'ematiques\\
Universit\'e de Gen\`eve\\
2--4 rue du Li\`evre\\
1211 Gen\`eve 24\\
Switzerland\\
\rm David.Cimasoni@math.unige.ch}
\vskip0.5cm

\nonumsection{References}

\end{document}